# THE ALPHA GROUP TENSORIAL METRIC

Cleber Souza Correa

Thiago Braido Nogueira de Melo

Diogo Machado Custódio
*Instituto de Aeronáutica e Espaço – IAE – Brasil*

## Abstract

The Alpha Group is an abstract geometry group in $R^4$. The way it was conceived allows a new interpretation of the structure of hypercomplex space with a new geometry and spatial topology, and a meaning for the geometric representation of $R^4$ space to infinity. Therefore, it has been described as the tensorial metric formula in the Alpha Group. It will be shown that the Riemannian and Euclidean distance metrics between infinitesimal surfaces are represented as special cases of the metric of the Alpha group.

**Keywords:** Abstract Algebras, Group Theory, Abstract Geometry.

## [MÉTRICA TENSORIAL NO GRUPO ALPHA]

### Resumo

O Grupo Alpha é um grupo de geometria abstrata no $R^4$. No qual propõe modificar, e dar uma nova interpretação na estrutura do espaço hiper-complexo. Criando uma nova estrutura geometria e de topologia espacial, como também, permite dar um significado na representação geométrica ao infinito no espaço do $R^4$. Portanto, foi apresentado neste artigo a descrição de uma fórmula métrica tensorial no Grupo Alpha. Busca-se mostrar que a métrica da distância entre superfícies infinitesimais de Riemanniana e Euclidiana são casos específicos da métrica do Grupo Alpha.

**Palavras-chave:** Álgebra Abstrata, Teoria de Grupo, Geometria Abstrata.

## Introduction - Historic View

The work of Georg Cantor (1884) and (1915) showed that infinite sets do not all have the same size and distinguished between enumerable and uncountable sets. Cantor proved that the set of rational numbers q is numerous (with cardinality equal to the set of natural numbers n). The set of real numbers R is uncountable (with a cardinality denoted by C greater than N), which allows the definition of a concept of different kinds of infinity. With the theoretical development and advances in mathematical knowledge in the 20th century, new mathematical theories were developed, such as the catastrophe theory of Thom, René (1972), which is based on the soft mapping theory of Whitney (Singularities (1955)) and the theory of bifurcations of dynamical systems of Poincaré, H. (1879) and Andronov, A. A. (1933).

The development of the ideas of abstract nature and mathematical thinking allowed the development of non-Euclidean geometries that broke with the paradigms of the postulates of Euclidean geometry. Lobachevsky's geometry showed that Euclid's geometry was not the absolute truth it was thought to be. In a sense, the discovery of non-Euclidean geometry dealt a crushing blow to Kantian philosophy, comparable to the effect it had on Pythagorean notions of the discovery of size. As a result of Lobachevsky's work, it became necessary to reexamine the fundamental ideas about the nature of mathematics (Halsted (1895)). Non-Euclidean geometry remained a rather marginal mathematics for some time until it was fully integrated into non-Euclidean geometry by the remarkably general ideas of G.F.B.Riemann (1826-1866) (Riemann and Weyl, 1923). Much more general than Lobachevsky's, where the question is simply how many parallels to a point are possible. Riemann recognized that geometry should not necessarily have to do with points, nor with straight lines, nor with space in the ordinary sense, but N-Uplas collections combined with certain rules also suggested a global view of geometry as the study of varieties of any number of dimensions in any kind of space. Riemann's work was of fundamental importance for geometrical ideas, his proposal for the general study of curved metric spaces and not for the special case of geometry over the sphere, which later became possible in general relativity, Boyer, Carl B. (1996).

In the field of abstract algebra, the studies of Sophus Lie (1842-1899), a contemporary student of Felix Klein (1849-1925) in Goettingen, wrote a three-volume treatise on the theory of transformation groups, Klein (1888 and 2004). Lie's contact transformations, systematized by Klein, established a bi-univocal correspondence between straight lines and spheres of Euclidean space, such that competing straight lines corresponded to tangent spheres. Such concepts, straight lines and spheres in three-dimensional Euclidean space form a four-dimensional space. Contact transformations are analytic transformations that transfer tangent surfaces to tangent surfaces. A collection of elements is said to form a group with respect to a particular operation if (1) the collection is closed under the operation, (2) the collection contains an identity element with respect to the operation, (3) each element in the collection contains an inverse element with respect to the operation, and (4) the operation is associative. The elements may be numbers (as in arithmetic), points (in geometry), transformations (in algebra or geometry), or something else. The operation can be arithmetic (like addition, multiplication, or division) or geometric (like a rotation about a point or axis). Or any other rule to combine two elements of a set (such as two transformations) to form a third element of the set. Here the generality of the group concept becomes clear. Therefore, this work attempts to characterize and define a tensor metric in the Alpha Group described by Correa, *et al.* (2022).

**The Alpha Group's Aspects**

According to group theory, the ratio of a ring (division) and a geometric (such as rotation about a point or axis), where there is also a rotation of one plane about the other 90 degrees, can be interpreted as a function of $\pi/2 + n\pi$. In which the rule combines two elements of a set, such as two transformations to form a third element of the set. This transformation is the operation performed to divide a multiplane by another with the same similarity. For all elements and the operation of a root with an index equal to 2. All the resulting operations are defined as a zero multiplication gives zero. As a result, this new structure appears as a region defining a relation between one plane and another, and a principal axis with multiples of imaginary μ is created. This operation can be interpreted as a ground for a canonical vector of the space of the Alpha group. It will represent the maximal deformation. It defines properties in this new group and creates a geometrical and topological structure. It will represent the maximum deformation and define properties in this new group. For the number μ (1/0) can be interpreted as a canonical vector. It is associated with the maximal deformation in $R^4$ defined by group theory as the ground for the quotient ring (Cohn (1982) and Jacobson (2009)) of two infinite complex planes. Such a topological and geometrical structure was described by Correa, *et al.* (2022) based on group theory.

**Results**

The Alpha Group is defined by its numerical structure by AG: a + bi + cμ + diμ, a real term (a), a second complex term with an imaginary number i (bi), the third term with the imaginary Alpha number μ(cμ) and the fourth term with two imaginary numbers iμ(diμ). The geometry of the Alpha group has a surface that you can see in Figure 1. It has several properties such as asymmetry and reflection of infinite numerical surfaces in $R^4$, where R and C represent the real and complex parts of the numerical structure of the alpha group. The imaginary number μ would represent the most fundamental relation of geometric infinity; all surfaces would be multiples of it.

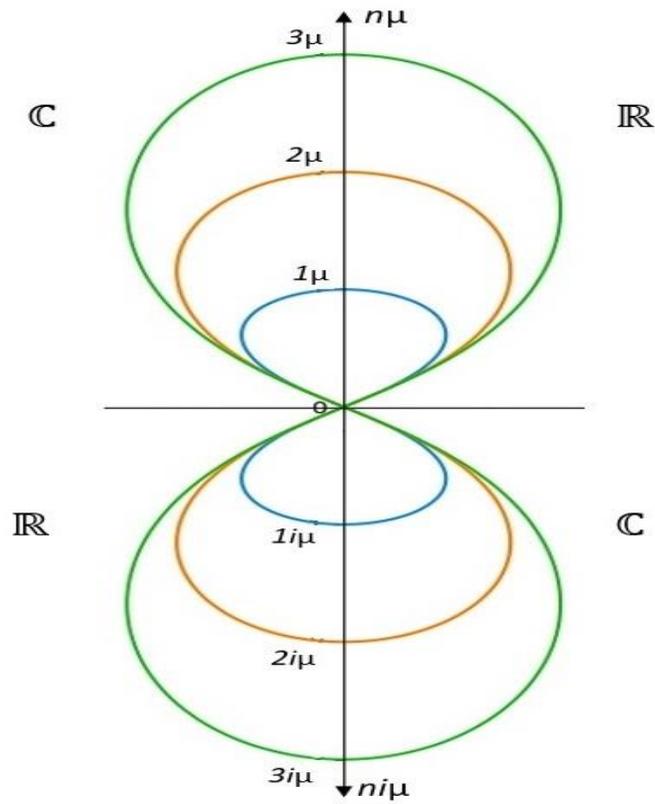

**Figure 1** - The geometric space of the Alpha Group in $R^4$, Poincaré cut.

In the case of a tensorial metric between two surfaces near the Alpha Group number of the $R^4$ group, it can be defined by a generalization of the distance formula of two infinitesimal points near the Riemann space in Equation I:

$ds^2 = g_{11} dx^2 + g_{12} dxdy \, i + g_{13} dxdz \, \mu + g_{14} dxdt \, i \, \mu + g_{21} dydx \, i + g_{22} dy^2 \, i^2 + g_{23} dydz \, i \, \mu + g_{24} dydt \, i^2 \, \mu \ + g_{31} dzdx \, \mu + g_{32} dzdy \, i \, \mu + g_{33} dz^2 \, \mu^2 + g_{34} dzdt \, i \, \mu^2 \ + g_{41} dtdx \, i \, \mu + g_{42} dtdy \, i^2\mu + g_{43} dtdz \, i \, \mu^2 + g_{44} dt^2 \, i^2\mu^2$ $\hspace{2cm}(I)$

With 16 terms in which gs can be constant or, more generally, functions of x, y, z, and t. Knowing that $i^2 = -1$ and $\mu^2 = \mu$ (the number $\mu$ has an infinity representation and the square is itself infinite). Equation II expresses the following:

$$ds^2 = g_{11} \ dx^2 + g_{12} \ dxdy \ i + g_{13} \ dxdz \ \mu + g_{14} \ dxdt \ i \ \mu$$
$$+ g_{21} \ dydx \ i - g_{22} \ dy^2 + g_{23} \ dydz \ i \ \mu - g_{24} \ dydt \ \mu$$
$$+ g_{31} \ dzdx \ \mu + g_{32} \ dzdy \ i \ \mu + g_{33} \ dz^2 \ \mu + g_{34} \ dzdt \ i \ \mu$$
$$+ g_{41} \ dtdx \ i \ \mu - g_{42} \ dtdy \ \mu + g_{43} \ dtdz \ i \ \mu - g_{44} \ dt^2 \ \mu$$
$$(II)$$

Regrouping the terms, we have Equation III:

$$ds^2 = g_{11} \ dx^2 - g_{22} \ dy^2$$
$$+ (g_{12} + g_{21}) \ dxdy \ i$$
$$+ (g_{13} \ dxdz - g_{24} \ dydt + g_{31} \ dzdx + g_{33} \ dz^2 - g_{42} \ dtdy - g_{44} \ dt^2) \ \mu$$
$$+ (g_{14} \ dtdx + g_{23} \ dydz + g_{32} \ dzdy + g_{34} \ dzdt + g_{41} \ dtdx + g_{43} \ dtdz) \ i \ \mu \quad (III)$$

This is the formula of tensorial metric between two surfaces in the Alpha Group. This would find the distance between two surfaces by a geodesic line connecting the two surfaces using all the terms in the equation III. This formula of the tensorial metric of the Alpha group has Riemannian space as a special case. Since we know that i and $\mu$ are constants, we set the terms $g_{14}$, $g_{24}$, $g_{34}$, $g_{41}$, $g_{42}$, $g_{43}$, and $g_{44}$ equal to zero within their respective gs. Since $g_{22}$ is considered a constant and we can exchange the negative signal for a positive signal, we get the equation IV:

$$ds^2 = g_{11} \ dx^2 + g_{12} \ dxdy + g_{13} \ dxdz$$
$$+ g_{21} \ dydx + g_{22} \ dy^2 + g_{23} \ dydz$$
$$+ g_{31} \ dzdx + g_{32} \ dzdy + g_{33} \ dz^2 \quad (IV)$$

*Regrouping the terms, we have equation V:*

$$ds^2 = g_{11} \ dx^2 + (g_{12} + g_{21}) \ dxdy + (g_{13} + g_{31}) \ dxdz$$
$$+ g_{22} \ dy^2 + (g_{23} + g_{32}) \ dydz$$
$$+ g_{33} \ dz^2 \quad (V)$$

Equation V represents a Riemannian space, is a specific case of the Alpha Group, continuing and considering $g_{11} = g_{22} = g_{33} = 1$ and other *gs* equal to *zero*.
We obtain equation VI characterizes Euclidean space.

$$ds^2 = g_{11} \ dx^2 + g_{22} \ dy^2 + g_{33} \ dz^2 \quad (VI)$$

**Conclusion**

In this paper, an attempt was made to show that the metrics of infinitesimal distance from Riemannian and Euclidean space are special cases of the metrics of the Alpha group. The Alpha group satisfies the properties defined by group theory. These results create a consistent structure that geometrically and topologically characterizes a hypercomplex numerical group. It has inherent properties that characterize transformations between surfaces. These definitions allow the development of new research in mathematics. The geometry of the alpha group is based on the theories of George Cantor, who was concerned about the nature of different kinds of infinity. This allows us to interpret the existence of the numerous and different types of infinity and their connection to geometry and topology.

**Cleber Souza Correa**
Instituto de Aeronáutica e Espaço

**E-mail:** clebercsc@fab.mil.br

**Thiago Braido Nogueira de Melo**
Instituto de Aeronáutica e Espaço

**E-mail:** braidotbnm@fab.mil.br

**Diogo Machado Custódio**
Instituto de Aeronáutica e Espaço

**E-mail:** diogodmc@fab.mil.br